Application of Renormalization Techniques to the Classical Arbelos Problem

by

Yelena Shvets, yelena.shvets@maine.edu

Jacob LiBrizzi, jal22@st-andrews.ac.uk

Abstract

This application of renormalization techniques offers a modern take on the classical Arbelos geometry problem. Keeping within the context of the original problem, two semicircles, meeting at chord T, are together circumscribed by a third semicircle. Separate from the original Arbelos result, both circumscribed semicircle areas are found in terms of chord T and the third circumscribing semicircle radius R. This approach eliminates the additional variables of the circumscribed semicircle radii.

Introduction

The following discussion expands on the classical Arbelos problem. The project began as a consequence of misreading, or reading beyond, a geometry final exam question. In the spring of 2016, at the University of Southern Maine, the famous Arbelos proposition appeared on the final exam given by Professor Silvia Valdes-Leon in her College Geometry course. The Arbelos, or "Shoemaker's Knife," proposition is accredited to the third century B.C. Greek mathematician Archimedes of Syracuse, and can be found as Proposition 4 in the *Book of Lemmas* (Archimedes 304). As seen in *The Works of Archimedes*, the proposition is stated as follows:

> If AB be the diameter of a semicircle and N any point on AB, and if semicircles be described within the first semicircle and having AN, BN as diameters respectively, the figure included between the circumferences of the three semicircles is 'what Archimedes called an ἄρβυλος';[1] and its area is equal to the circle on PN as diameter, where PN is perpendicular to AB and meets the original semicircle in P. (304)

Defining *PN*'s length as *T*, the exam asked for the area bounded between all three semicircles. This region has been known as the Shoemaker's Knife. As shown below, letting θ represent the Shoemaker's Knife, the area of θ is:

$$A(\theta) = \frac{\pi T^2}{4} = \pi \left(\frac{T}{2}\right)^2.$$

---

[1] The name Arbelos comes from the Greek ἄρβυλος meaning Shoemaker's Knife.

## Traditional Arbelos Proof

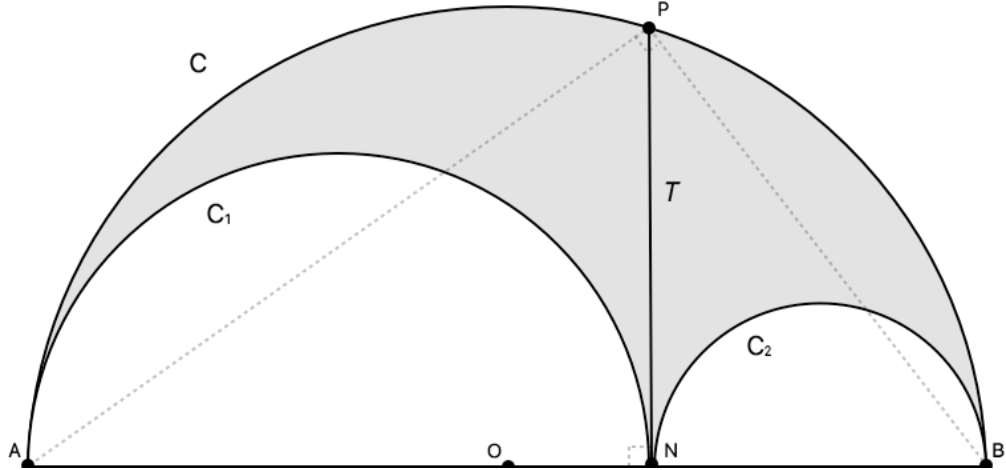

1. Let $R, R_1, R_2$ be radii of semicircles $C, C_1, C_2$ respectively.

   Then $R = R_1 + R_2$ $\Rightarrow$ $R^2 = (R_1 + R_2)^2$ $\Rightarrow$ $(R_1)^2 + (R_2)^2 = R^2 - 2R_1R_2$.

   We have: i) $AN = 2R_1$, $NB = 2R_2$, $PN = T$, $m(\angle ANP) = 90°$

   ii) By properties of circles $\angle APB = \frac{1}{2}\angle AOB$, (where $AO$ is a radius of $C$)

   $\Rightarrow$ Triangle $\triangle APB$ is a right triangle at $P$.

2. Let the area of a semicircle $C$ be denoted by $A(C)$.

   Then $A(C) = A(C_1) + A(C_2) + A(\theta)$,     (where $A(\theta)$ denotes the area of the
   
   Shoemaker's Knife region bound
   
   between semicircles $C, C_1, C_2$)

   and $A(\theta) = A(C) - \big(A(C_1) + A(C_2)\big)$.

3. $(AB)^2 = (AP)^2 + (PB)^2 = [(AN)^2 + (PN)^2] + [(NB)^2 + (PN)^2]$

$$= (AN)^2 + (NB)^2 + 2(PN)^2$$

and $(AB)^2 = (AN + NB)^2 = [(AN)^2 + (NB)^2 + 2(AN)(NB)]$.

Then $2(PN)^2 = 2(AN)(NB) \Rightarrow (PN)^2 = (AN)(NB)$

$$\Rightarrow T^2 = (2R_1)(2R_2) = 4R_1R_2 \; (by \; 1.)$$

$$\Rightarrow R_1R_2 = \frac{T^2}{4}.$$

4. Then $A(\theta) = \frac{1}{2}\pi R^2 - \left(\frac{1}{2}\pi(R_1)^2 + \frac{1}{2}\pi(R_2)^2\right)$

$$= \frac{\pi}{2}[R^2 - ((R_1)^2 + (R_2)^2)]$$

$$= \frac{\pi}{2}[R^2 - (R^2 - 2R_1R_2)]$$

$$= \frac{\pi}{2}[2R_1R_2]$$

$$= \pi R_1R_2$$

$$= \frac{\pi}{4}T^2 \; (by \; Part \; 3).$$

∴ The area of Shoemaker's Knife region is $A(\theta) = \pi\left(\frac{T}{2}\right)^2$ ∎

This result is well known, if somewhat surprising, in that it appears, on the surface, that the area of the Arbelos does not depend on the radius of the circumscribing semicircle $C$. This is of course, misleading, since the length of chord $T$ is bound the radius: $0 \leq T \leq R$.

The misreading, or reading beyond, occurred when the second author assumed that he also needed to find the areas of each of the circumscribed semicircles in terms of the variable $T$ and the radius R of the circumscribing semicircle. This turned out to be an elusively simple matter, upon the introduction of a very natural renormalization, or rescaling (Barenblatt Ch. 1: section 1.1 and Ch. 5: section 5.3).

<p align="center">Renormalization and Variable Reduction</p>

It's easy to see that the sum of the areas of the two semicircles depends on the radius of the circumscribing semicircle and the chord $T$:

$$A(C_1) + A(C_2) = A(C) - A(\theta) = \frac{\pi}{2}R^2 - \frac{\pi}{4}T^2.$$

In addition, the sum of the two smaller radii is equal to the radius R.

We have now an opportunity to reduce the number of variables:

$$R = R_1 + R_2$$

$$\frac{R}{R} = \frac{R_1}{R} + \frac{R_2}{R}$$

$$1 = \frac{R_1}{R} + \frac{R_2}{R}$$

Letting $r_1 = \frac{R_1}{R}$ and $r_2 = \frac{R_2}{R}$, we have that $r_2 + r_1 = 1$ and $r_2 = 1 - r_1$.

We also introduce a variable $t = \frac{T}{R}$, obtaining a set of two dimensionless parameters:

$r_1$ and $t$, which can now be connected using (3) above:

$$T^2 = 4R_1R_2$$

$$\frac{T^2}{R^2} = \frac{4R_1R_2}{R^2}$$

$$t^2 = 4r_1r_2$$

$$t^2 = 4r_1(1-r_1)$$

$$\frac{t^2}{4} = r_1 - r_1^2$$

$$\frac{1-t^2}{4} = \left(r_1 - \frac{1}{2}\right)^2$$

We can solve this quadratic for $r_1$, using t as a parameter, and obtain:

$$r_1 = \frac{1 \pm \sqrt{1-t^2}}{2}$$

If we select $r_1 = \frac{1+\sqrt{1-t^2}}{2}$, then $r_2 = \frac{1-\sqrt{1-t^2}}{2}$, and vice-versa. In this way, each of the smaller radii is uniquely determined by $t$, or by the ratio of $T$ to $R$.

Clearly, the area of each of the semicircles $C_1$ and $C_2$ can now be expressed as follows:

$A(C_1) = \frac{\pi}{2}(R_1)^2 = \frac{\pi}{2}(Rr_1)^2 = \pi R^2 \frac{\left(1+\sqrt{1-t^2}\right)^2}{8} = \frac{\pi}{8}\left(R + \sqrt{R^2 - T^2}\right)^2$; and similarly

$A(C_2) = \frac{\pi}{8}\left(R - \sqrt{R^2 - T^2}\right)^2$.

It's important to notice that $T$ varies between 0 and $R$, so that the expressions above are always defined.

Acknowledgments

The authors are grateful to Professor Silvia Valdes-Leon who introduced them to the Arbelos Problem and encourage them to write this note.